\documentclass[11pt,reqno]{amsart}

\usepackage[T1]{fontenc}
\usepackage{mathptmx,microtype,hyperref,graphicx,dsfont}
\usepackage{amsthm,amsmath,amssymb}
\usepackage[foot]{amsaddr}
\usepackage{enumitem}
\usepackage[nameinlink,capitalise]{cleveref}
\usepackage[font=footnotesize,margin=2cm]{caption}

\crefname{theorem}{Theorem}{Theorems}
\crefname{thm}{Theorem}{Theorems}
\crefname{mainthm}{Theorem}{Theorems}
\crefname{lemma}{Lemma}{Lemmas}
\crefname{lem}{Lemma}{Lemmas}
\crefname{remark}{Remark}{Remarks}
\crefname{claim}{Claim}{Claims}
\crefname{subclaim}{Sub-claim}{Sub-claims}
\crefname{prop}{Proposition}{Propositions}
\crefname{proposition}{Proposition}{Propositions}
\crefname{defn}{Definition}{Definitions}
\crefname{corollary}{Corollary}{Corollaries}
\crefname{conjecture}{Conjecture}{Conjectures}
\crefname{question}{Question}{Questions}
\crefname{chapter}{Chapter}{Chapters}
\crefname{section}{Section}{Sections}
\crefname{figure}{Figure}{Figures}
\newcommand{\crefrangeconjunction}{--}

\theoremstyle{plain}

\newtheorem{thm}{Theorem}
\newtheorem*{thm*}{Theorem}

\theoremstyle{definition}

\theoremstyle{remark}

\renewcommand{\P}{{\bf P}}
\newcommand{\E}{{\bf E}}

\newcommand{\Z}{{\mathbb Z}}

\newcommand{\G}{{\mathcal G}}
\newcommand{\Gnp}{\G_{n,p}}

\author[B. Kolesnik]{Brett Kolesnik}

\address{Department of Statistics, University of California, Berkeley}

\email{bkolesnik@berkeley.edu}

\begin{document}

\title[Large deviations of the greedy algorithm]
{
Large deviations of the greedy independent set algorithm 
on sparse random graphs 
}

\begin{abstract}
We study the greedy
independent set algorithm on sparse 
Erd{\H{o}}s--R{\'e}nyi random graphs $\G(n,c/n)$. 
This range of $p$ is of interest due to the 
threshold at $c=e$, beyond which 
it appears that greedy algorithms 
are affected by a sudden change in the 
independent set landscape. 
A large deviation principle was recently established by Bermolen et al. (2020), 
however, the proof and rate function are somewhat involved. 
Upper bounds for the rate function were obtained earlier by Pittel (1982). 
By discrete calculus, 
we identify the optimal 
trajectory realizing a given large deviation
and obtain the rate function in a simple closed form. 
In particular, we show that Pittel's bounds are sharp. 
The proof is brief and elementary. 
We think the methods presented here will be useful 
in analyzing  the tail behavior of 
other random growth and exploration processes. 
\end{abstract}

\maketitle

\section{Introduction}\label{S_intro}

We investigate the size of the independent set 
found by the greedy algorithm 
on the sparse Erd{\H{o}}s--R{\'e}nyi \cite{ER59} graph $\G(n,c/n)$. 
Each pair of vertices in $\G$ is joined by an edge
independently with probability $p=c/n$, 
where $c\in(0,\infty)$ is a constant, fixed throughout this work. 
Recall that a set $I\subset [n]$ is
{\it independent}  in $\G$ if no two vertices in $I$
are neighbors. 

The {\it greedy algorithm} is a local exploration process.
Initially, all vertices are {\it available,}  $A_0=[n]$. 
Beginning with $I_0=\emptyset$, 
in each step $k\ge1$, an independent set $I_k$ is formed by
adding to $I_{k-1}$ a random vertex 
$v_k$ in the set $A_{k-1}$ of currently available 
vertices with no neighbors in $I_{k-1}$. Edges from $v_k$
are then revealed, and neighbors of $v_k$ in $A_{k-1}$ are removed 
to obtain $A_{k}$. 
The size of the independent set eventually obtained
is the first step $S_g=\min\{k\ge0:A_k=\emptyset\}$ 
in which  $I_k$ is maximal. 

As shown by Pittel~\cite{P82}, a set of size approximately $s_cn$ is 
typically found, where $s_c=(1/c)\log(1+c)$. 
In this work, we focus on the atypical behavior. 
For $s\neq s_c$, we let $P_s$ denote $\P(S_g\le sn)$ if $s<s_c$
and $\P(S_g\ge sn)$ if $s>s_c$. 

\begin{thm}\label{T_main}
Fix $s\neq s_c$. Suppose that $s_n\to s$ 
as $n\to\infty$. 
Then 
\begin{equation}\label{E_rate}
\lim_{n\to\infty}\frac{1}{n}\log P_{s_n}= 
\log a_s
+\frac{1}{c}\int_{a_s}^{b_s}\frac{\log u}{1-u}du
\end{equation}
where $b_s=a_se^{c(1-1/a_s)}$ and $a_s>0$ uniquely satisfies 
\begin{equation}\label{E_a}
s=\frac{1}{c}\log \frac{b_s-1}{a_s-1}.
\end{equation}
\end{thm}

The right hand side of \eqref{E_a}
converges to $s_c$  
as $a\to1$. 
Otherwise, $a\in(0,1)$ if $s<s_c$ and $a\in(1,\infty)$
if $s>s_c$. 

The upper bounds in \cref{T_main} 
are obtained in \cite{P82} by martingale arguments. 
Indeed, setting $\delta=c$, $x=s$ and $y=a_s$ in Lemma 4 of \cite{P82}  
yields the upper bound $F_c(s,a_s)$ (in the notation of \cite{P82}), 
equal to the right hand side in \eqref{E_rate}. 
Our proof of \cref{T_main}, by discrete calculus, follows a completely 
different approach. Futhermore, \cref{T_main} 
shows that these upper bounds are sharp.
We think the methods of the current article 
are of independent interest, and 
will be useful 
in analyzing  
a variety of other random growth and exploration processes. 

\cref{T_main} 
is obtained in part by identifying the least-cost trajectory 
realizing a given deviation from $s_c$. 
The key ingredient in the proof is a 
discrete analogue of the Euler--Lagrange equation
due to Guseinov \cite{G04}. 

\begin{thm}\label{T_opt}
Fix $s\neq s_c$. Asymptotically, 
\[
\hat y_s(x)
=\frac{1+(a_s-1)e^{cx}}{e^{cx}}
[\frac{1}{a_s}-\frac{1}{c(a_s-1)}\log\frac{1+(a_s-1)e^{cx}}{a_s}]
\]
is the optimal trajectory 
for 
$|A_{xn}|/n$ 
amongst those
decreasing from 1 to 0 over $[0,s]$. 
More specifically, if $s_n\to s$, then 
\[
\lim_{n\to\infty}\frac{1}{n}\log P_{s_n}= \int_0^s \ell_s(x)dx
\]
where
\begin{align*}
\ell_s(x)
&=-(1+\hat y'_s(x))[1+\frac{c\hat y_s(x)}{1+\hat y'_s(x)}
+\log(-\frac{c\hat y_s(x)}{1+\hat y'_s(x)})]\\
&=
[\log\frac{1+(a_s-1)e^{cx}}{a_s}-\frac{c(a_s-1)}{a_s}]
[1-\frac{\log(1+(a_s-1)e^{cx})}{(a_s-1)e^{cx}}]
\end{align*}
is the cost function associated with 
$\hat y_s$. 
\end{thm}

\begin{figure}[h]
\centering
\includegraphics[scale=0.26]{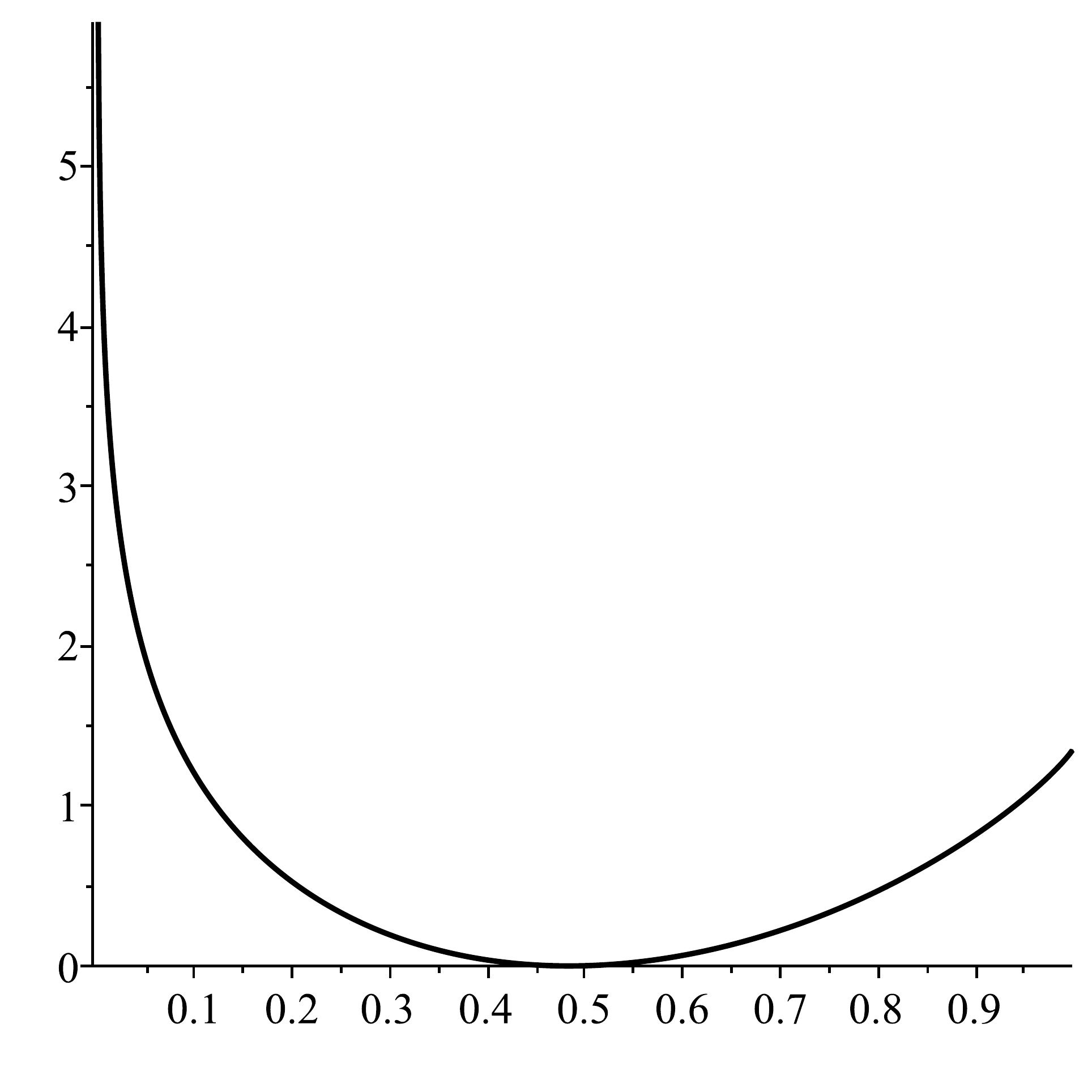}
\includegraphics[scale=0.26]{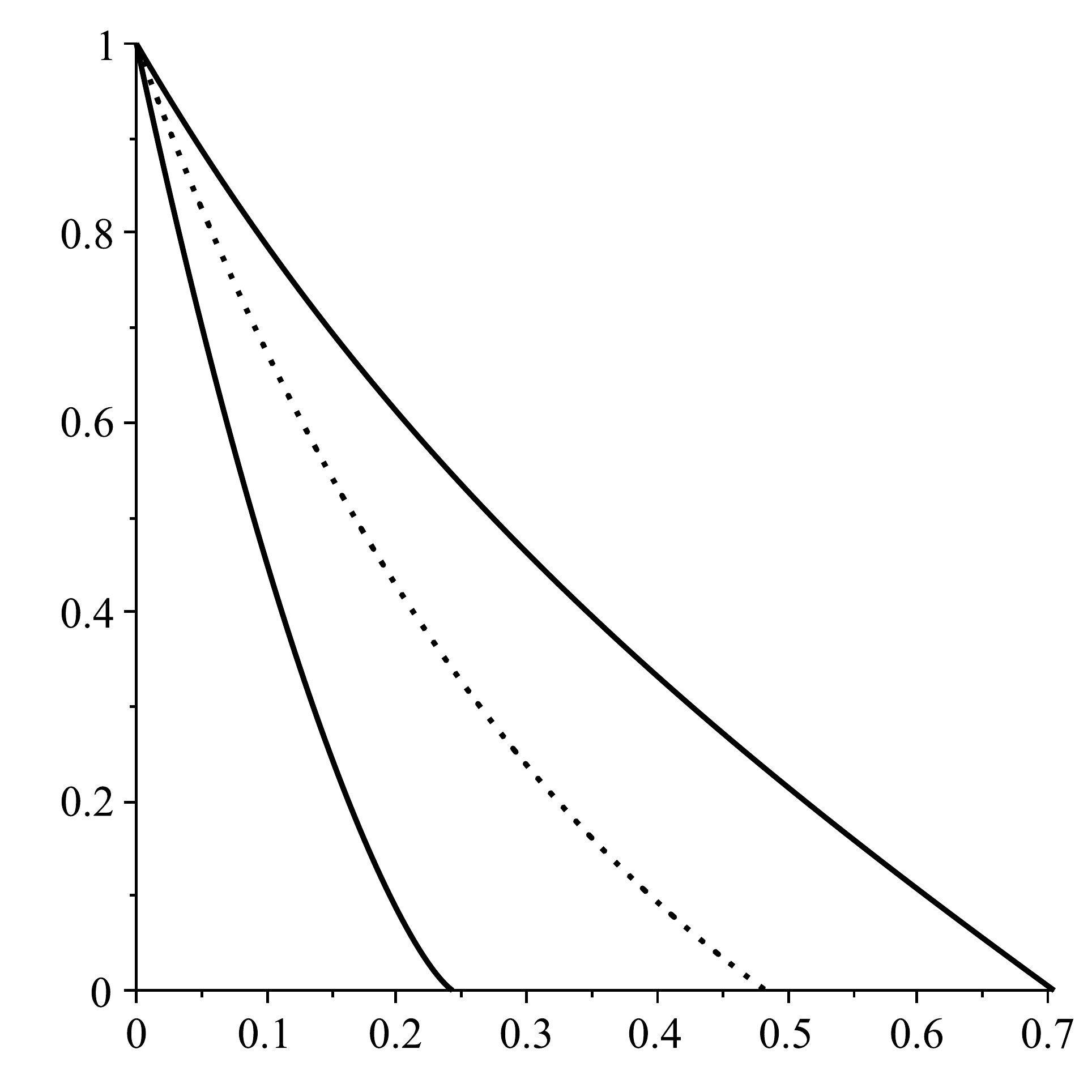}
\vspace{-0.25cm}
\caption{
In both figures, $c=e$. 
The rate function is at left, 
intersecting the $s$-axis at 
$s_c\approx 0.483$. 
At right, the expected trajectory $\bar y_{c}$ 
is dotted, between two least-cost  
deviating trajectories
$\hat y_s$, associated with 
the values $a=1/2$ ($s\approx 0.243$) and $a=2$ ($s\approx 0.704$).
}
\label{F_lam}
\end{figure}

Taking $a\to1$, we recover 
the typical  
trajectory
\[
\bar y_c(x)=
(1+1/c)e^{-cx}-1/c
\]
which decreases from 1 to 0 over $[0,s_c]$.

Note that $\ell_s(x)=-\Gamma^*_\lambda(\xi)$, 
where 
$\Gamma^*_\lambda(\xi)
=\sup_\vartheta[\vartheta \xi-\Gamma_\lambda(\vartheta)]$
is the Legendre--Fenchel transformation 
of the cumulant-generating function $\Gamma_\lambda(\vartheta)=\lambda(e^\vartheta-1)$
of a rate $\lambda$ Poisson random 
variable, evaluated at 
$\lambda_x=c\hat y_s(x)$ and
$\xi_x=-(1+\hat y'_s(x))$. 

Detailed heuristics for $s_c$, $\bar y_{c}$, $\hat y_s$
and $\ell_s$ are given in \cref{S_heuristics} below.

\subsection{Discussion}

Recent work by Bermolen et al.\ \cite{BGJM20} proves
a large deviation principle for the scaled trajectory $y_s$  of the
greedy algorithm and the size $S_g$ of the
independent set obtained. This is based
on the approach of Feng and Kurtz \cite{FK06} to 
large deviations by  
the theory of viscosity solutions and 
convergence of nonlinear semigroups. 
The arguments in \cite{BGJM20}  
are somewhat lengthy and complicated
and the rate function given there, although explicit, 
is not as convenient as possible. 

We used Guseinov's \cite{G04}
discretized Euler--Lagrange equation
recently \cite{AK18} to derive sharp tail estimates for a
certain percolation model.  
We suspect that the utility of  \cite{G04} 
is wide-ranging, although it is 
perhaps not well-known
in the context of discrete probability. 
In the current context, it leads to a natural proof: 
The heuristic in \cref{S_opt} 
correctly identifies the optimal trajectory $\hat y_s$, using Poisson approximation
to the Binomial and 
the usual Euler--Lagrange equation. The discretized version 
\cite{G04} provides the means to make this rigorous (\cref{S_proof}).

\subsection{Motivation}

The 
range $p=\Theta(1/n)$ is of interest 
in relation to greedy algorithms 
due to the so-called
{\it $e$-cutoff phenomenon}. 
As discussed in \cite{BGJM20}, 
one might hope 
\cite{T09}
that further analysis of the large deviations
of such algorithms could bring 
this threshold into clearer light. 
See \cite{BGJM20} for some loose observations in this direction.
It is known \cite{KS81,AFP98} that a slight modification of the greedy algorithm, 
called the {\it degree-greedy algorithm,}
almost surely finds an optimal independent set for $c< e$. 
However, when $c>e$, the situation is more complicated. 
For instance, showing  
that the scaled expected independence number
$\E({\mathcal I}_{c,n}/n)$ of $\G(n,c/n)$ has a limiting value
remained open for some time (see e.g.\ \cite{AldousOP}), until the breakthrough \cite{BGT13}. 
See also \cite{GNS06}, where a more  amenable
weighted version is analyzed.

\subsection{Acknowledgements}
We thank Boris Pittel for useful discussions.
We thank 
Shirshendu Ganguly for bringing the
graphical description of the greedy algorithm to our attention, 
during conversations together with Daniel Reichman. 

\section{Heuristics}\label{S_heuristics}

\subsection{Typical behavior}\label{S_typ}
A simple heuristic for  
$s_c$ and $\bar y_{c}$ can be seen  
by 
a graphical description of the greedy algorithm. 
Although the proof in \cref{S_proof} does not rely on this description, 
it lends some useful intuition. 
See also McDiarmid~\cite{McD84} for an alternative explanation.
Similar arguments are used by Nachmias and Peres \cite{NP10}, for instance, 
in analyzing the size of the largest component in 
the critical graph $\G(n,1/n)$. This analysis
is based on the exploration processes of 
Karp \cite{K90} and Martin-L\"{o}f \cite{ML86}. 

Start with a row of $n$
particles, one for each vertex in $\Gnp$. In step $k$, 
and unmarked particle is marked, and then all 
particles move up to the next row independently 
with probability $1-p$. 
The independent 
set obtained by the greedy algorithm consists of all 
marked particles, its size $S_g$ being  
the index of the first row without any unmarked particles. The 
available vertices $A_k$ are those in the $k$th
row that are unmarked.

In the $k$th row, we expect $n(1-p)^k$ particles, and 
$\sum_{\ell=1}^k (1-p)^\ell$ of these to be marked. Therefore, 
we expect 
\[
n(1-p)^k-(1/p-1)[1-(1-p)^k]
\sim n \bar y_{c}(k/n)
\] 
unmarked particles in the $k$th row. 
This diverges to either $+\infty$ or $-\infty$, 
depending on whether $k/n<s_c$ or $k/n>s_c$. 
Hence it can be shown that $\E (S_g/n)\sim s_c$.
This was first proved by Pittel \cite{P82} (cf.\ \cite{GMcD75,K76,McD84,BKS17}).  

\subsection{Large deviations}\label{S_opt}
The optimal deviating trajectory $\hat y_s(x)$
given in \cref{T_opt} can be guessed by Poisson approximation
and the  Euler--Lagrange equation. 

Consider some trajectory $y_s(x)$ for $|A_{xn}|/n$, decreasing 
from $1$ to 0 over $[0,s]$, leading to an independent set of size
$sn$. 
Suppose that $|A_{xn}|$ has followed this trajectory up until $x=k/n$.
Before proceeding to the next step, 
one of the $ny_s(x)$ available particles $v_k$ is marked (see \cref{S_typ} above). 
Approximately a Poisson
with rate $c y_s(x)=\lambda_x$ of the rest
are neighbors with $v_k$. Hence, to continue along this trajectory until
$x'=x+1/n$, we require this random variable to take the 
value $ny_s(x)-1-ny_s(x')\approx -(1+y'_s(x))=\xi_x$. 
The log probability of this event is approximately 
$-\Gamma^*_{\lambda_x}(\xi_x)=\ell_s(x)$. Hence
the log probability that $|A_{xn}|/n$ follows $y_s(x)$ over $[0,s]$
is approximately $n\int_0^s\ell_s(x)dx$. 
Therefore, by the Euler--Lagrange equation, 
we expect the optimal trajectory
$\hat y_s(x)$ to satisfy 
\[
\frac{1+\hat y'_s(x)}{\hat y_s(x)}+c
=\frac{d}{dx}\log(-\frac{c \hat y_s(x)}{1+ \hat y'_s(x)}). 
\] 
It can be seen that 
$\hat y_s(x)$ 
solves this equation subject to  
$y_s(0)=1$ and $y_s(s)=0$.

\section{The proof}\label{S_proof}

First we show that, in the limit, 
the tail probability $P_s$ is dominated by a single optimal trajectory
$\hat y_s$, which we identify by discrete calculus
of variations. 
This establishes the upper bound. The matching lower bound follows by 
considering any given trajectory that is sufficiently close to  
$\hat y_s$. 

\subsection{Upper bound}\label{S_UB}
Let $s_n\to s\neq s_c$ be given. 
Let us assume that $s>s_c$ and so $a_s>1$ in \eqref{E_a}. 
The same argument works for $s<s_c$,
with only minor changes.

For any $t\ge s_n$, 
let 
$0=x_0<x_1<\cdots<x_m=t$ 
be evenly spaced points (to the extent possible 
subject to all $x_in\in \Z$) where $m$ is chosen 
so that all $\Delta x_i=x_{i+1}-x_i=\Theta[(\log n)^2/n]$.
The choice of $(\log n)^2$ is not important, 
only that $m\log n\ll n$.

Let ${\mathcal Y}_{t}\subset\Z^{m+1}$ denote the set of possible values 
$n=Y_0>Y_1>\cdots>Y_m=0$ taken by the sequence 
of $|A_{x_i n}|$. 
All relevant trajectories are strictly decreasing, 
since at least one available vertex is removed in each step. 
Put ${\mathcal Y}_{s_n}^+=\bigcup_{t\ge s_n}{\mathcal Y}_{t}$. 
Then, taking a union bound,  
\[
P_{s_n}\le \sum_{Y\in{\mathcal Y}_{s_n}^+ }\prod_{i=0}^{m-1}
P_i(Y)
\]
where
\[
P_i(Y)=\P(|A_{x_{i+1} n}|=Y_{i+1}||A_{x_{i} n}|=Y_{i}).
\]
Moreover, since   
$|{\mathcal Y}_{s_n}^+|\le O(n^m)$, it follows,  
by the choice of $m$, that 
\begin{equation}\label{E_hatx}
\frac{1}{n}\log P_{s_n}\le 
o(1)+
\frac{1}{n}\sum_{i=0}^{m-1}
\log P_i(\hat Y)
\end{equation}
where $\hat Y$ maximizes 
$\sum_i
\log P_i(Y)$
over $Y\in{\mathcal Y}_{s_n}^+$. 

For all relevant $Y$, we have that 
\[
P_i(Y)
\le {Y_{i}-n\Delta x_i\choose Y_{i+1}}[(1-p)^{n\Delta x_i}]^{Y_{i+1}}
[1-(1-p)^{n\Delta x_i}]^{-(n\Delta x_i+\Delta Y_i)}. 
\]
Therefore, since $1-(1-p)^{n\delta}=c\delta(1+O(\delta))$, 
and using the standard bounds 
${k\choose \ell}\le(ek/(k-\ell))^{k-\ell}$
and $(1-x)^y\le e^{-xy}$, 
we find that 
\begin{equation}\label{E_PxiUB}
P_i(Y)
\le
[-e\frac{cY_{i}}{n+\Delta Y_{i}/\Delta x_i}
(1+O(\Delta x_i))]^{-(n\Delta x_i+\Delta Y_{i})}
e^{-cY_{i} \Delta x_i }.
\end{equation}
Observe that
\[
\frac{1}{n}\sum_{i=1}^m (n\Delta x_i+\Delta Y_{i})
\log(1+O(\Delta x_i))
\le O[\log(1+O(1/m))]\ll1. 
\]
Therefore 
\begin{equation}\label{E_f}
\frac{1}{n}\sum_{i=0}^{m-1}
\log P_i(Y)\le
o(1)
+\sum_{i=0}^{m-1} f(y_{i},\Delta y_{i}/\Delta x_i)\Delta x_i
\end{equation}
where $y=Y/n$ and 
\[
f(u,w)=-(1+w)[1+\frac{cu}{1+w}+\log(-\frac{cu}{1+w})].
\]
In upper bounding the righthand side of \eqref{E_f}, 
we may relax the restriction that all $y_i n\in\Z$,
and instead optimize over real-valued $y$.
Then, applying Theorem 5 in \cite{G04}, 
we find that the maximizer $\hat y$ satisfies 
\begin{equation}\label{E_Euler}
   f_u(\hat y_{i+1},\Delta \hat y_{i+1}/\Delta x_{i+1})
=
\Delta  f_w(\hat y_{i},\Delta \hat y_{i}/\Delta x_i)/\Delta x_i.
\end{equation}
Note that 
\[
f_u=-\frac{1+w}{u}-c,\quad 
f_w=-\log(-\frac{cu}{1+w}). 
\]
Standard results on Euler's method   
(see e.g.\ Section I.7 in \cite{HNW93},
in particular Theorems 7.3 and 7.5), 
imply that $\hat y_i$ and $\Delta \hat y_i/\Delta x_i$
are within $O(1/m)$ of $\hat y_t$
and $\hat y_t'$, where
$\hat y_t$ is the limiting trajectory 
satisfying 
\begin{equation}\label{E_DE}
\frac{1+\hat y'_t(x)}{\hat y_t(x)}+c=\frac{d}{dx}\log(-\frac{c \hat y_t(x)}{1+ \hat y'_t(x)})
\end{equation}
subject to $\hat y_t(0)=1$ and $\hat y_t(t)=0$. 
Therefore, by \eqref{E_hatx}, \eqref{E_f} and \eqref{E_Euler}, 
\begin{equation}\label{E_ubint}
\lim_{n\to\infty}\frac{1}{n}\log P_{s_n}
\le \max_{t\in[s,1]}\int_0^t f(\hat y_t(x),\hat y_t'(x))dx. 
\end{equation}

To solve \eqref{E_DE}, we first observe that 
\[
c(z(x)+1)=\frac{d}{dx}\log(-1/z(x))\implies z(x)=-\frac{1}{1+(a-1)e^{cx}}
\]
and then that
\begin{multline*}
1+\hat y'(x)=-\frac{c\hat y(x)}{1+(a-1)e^{cx}}\\
\implies
\hat y(x)=\frac{1+(a-1)e^{cx}}{e^{cx}}
[b-\frac{1}{c(a-1)}\log(1+(a-1)e^{cx})].
\end{multline*}
The boundary conditions imply that 
\[
a=a_t,\quad b=\frac{\log a_t}{c(a_t-1)}+\frac{1}{a_t}.
\]
Hence 
\begin{equation}\label{E_haty}
\hat y_t(x)
=\frac{1+(a_t-1)e^{cx}}{e^{cx}}
[\frac{1}{a_t}-\frac{1}{c(a_t-1)}\log\frac{1+(a_t-1)e^{cx}}{a_t}]
\end{equation}
as appears in \cref{T_opt}.

Next, in order to evaluate the integral in \eqref{E_ubint}, 
note that 
\[
-\frac{c\hat y_t(x)}{1+\hat y_t'(x)}=1+(a_t-1)e^{cx} 
\]
and
\[
1+\hat y_t'(x)
=
\frac{1}{(a_t-1)e^{cx}}[\log\frac{1+(a_t-1)e^{cx}}{a_t}-\frac{c(a_t-1)}{a_t}].
\]
Hence
\[
f(\hat y_t(x),\hat y_t'(x))=
[\log\frac{1+(a_t-1)e^{cx}}{a_t}-\frac{c(a_t-1)}{a_t}]
[1-\frac{\log(1+(a_t-1)e^{cx})}{(a_t-1)e^{cx}}].
\]
Note that this is $\ell_t(x)$ in  
\cref{T_opt}. 
Then, by basic calculus, 
\[
\int f(\hat y_t(x),\hat y_t'(x))dx
=-\hat y_t(x)\log(1+(a_t-1)e^{cx})+
\frac{1}{c}\int_1^{1+(a_t-1)e^{cx}}\frac{\log u}{1-u}du.
\]
Therefore, since $1+(a_t-1)e^{ct}=b_t$ by \eqref{E_a}, 
it follows
by \eqref{E_ubint} and \eqref{E_haty}
 that 
\begin{equation}\label{E_ubint2}
\lim_{n\to\infty}\frac{1}{n}\log P_{s_n}
\le \max_{t\in[s,1]}[\log a_t+
\frac{1}{c}\int_{a_t}^{b_t}\frac{\log u}{1-u}du].
\end{equation}

Finally, to complete the proof we show that 
$t=s$ is the maximizing case, 
as is expected, since $s$ is the least
extreme deviation from $s_c$ in $[s,1]$. 
Note that $t>s_c$
correspond with $a_t\in(1,\infty)$ in \eqref{E_a}, 
which are increasing in $t$. 
Hence we show that the right hand side 
in \eqref{E_ubint2} 
is decreasing in $a>1$. 
Observe that, with $b=ae^{c(1-1/a)}$,  
\[
\frac{d}{da}[\log a+
\frac{1}{c}\int_{a}^{b}\frac{\log u}{1-u}du]
=\frac{(1-c(1-1/a))(b/a)-1}{b-1}
(\frac{\log a}{c(a-1)}+\frac{1}{a}).
\]
Since 
\[
\frac{\log a}{c(a-1)}+\frac{1}{a}>0
\]
for all $a>0$, and $b>a>1$ for $a>1$, 
we need only check that 
\[
(1-\frac{c(a-1)}{a})\frac{b}{a}<1.
\] 
This follows noting that the left hand side is equal to 1 when 
$a=1$, since then also $b=1$,  and that 
\[
\frac{d}{da}[(1-\frac{c(a-1)}{a})\frac{b}{a}]=\frac{(1-a)bc^2}{a^4}<0.
\]
This proves the claim. 

Altogether, by \eqref{E_ubint2}, we find that 
\[
\lim_{n\to\infty}\frac{1}{n}\log P_{s_n}
\le \log a_s+
\frac{1}{c}\int_{a_s}^{b_s}\frac{\log u}{1-u}du
\]
as required. 

\subsection{Lower bound}\label{S_LB}
Having identified $\hat y_s$, 
the lower bound
follows easily by considering any trajectory of 
$|A_{xn}|/n$ sufficiently close to $\hat y_s$. 
We sketch the argument (using the notation from \cref{S_UB}).

The first step is to find a lower bound for $P_{i}(Y)$
which agrees with the upper bound \eqref{E_PxiUB}
up to small error. 
For all relevant $Y$, 
\begin{align*}
P_i(Y)
&\ge {Y_{i}-n\Delta x_i\choose Y_{i+1}}
(1-p)^{{n\Delta x_i\choose2}}
[(1-p)^{n\Delta x_i}]^{Y_{i+1}}
[np\Delta x_i (1-p)^{n\Delta x_i-1}]^{-(n\Delta x_i+\Delta Y_i)}\\
&\ge
{Y_{i}-n\Delta x_i\choose -(n\Delta x_i+\Delta Y_i)}
(1-p)^{nY_{i+1}\Delta x_i}
(np\Delta x_i)^{-(n\Delta x_i+\Delta Y_i)}. 
\end{align*}
Using the bound 
${k\choose \ell}\ge(e(k/\ell-1))^\ell/e\ell$, it follows that 
\[
P_i(Y)\ge 
(-e\frac{cY_{i+1}}{n+\Delta Y_{i}/\Delta x_i})^{-(n\Delta x_i+\Delta Y_{i})}
\frac{(1-p)^{nY_{i+1}\Delta x_i}}{-e(n\Delta x_i+\Delta Y_{i})}. 
\]
Therefore, by the choice of $m$, this implies that 
\[
\frac{1}{n}\sum_{i=0}^{m-1}
\log P_i(Y)\ge
o(1)
+\sum_{i=0}^{m-1} f(y_{i+1},\Delta y_{i}/\Delta x_i)\Delta x_i. 
\]
Hence, comparing this with \eqref{E_f}, it is not 
hard to see, by considering any possible trajectory 
$y$ sufficiently close to $\hat y_s$, that 
\[
\lim_{n\to\infty}\frac{1}{n}\log P_{s_n}
\ge
\log a_s+
\frac{1}{c}\int_{a_s}^{b_s}\frac{\log u}{1-u}du.
\]

\providecommand{\bysame}{\leavevmode\hbox to3em{\hrulefill}\thinspace}
\providecommand{\MR}{\relax\ifhmode\unskip\space\fi MR }
\providecommand{\MRhref}[2]{%
  \href{http://www.ams.org/mathscinet-getitem?mr=#1}{#2}
}
\providecommand{\href}[2]{#2}

\end{document}